\documentclass[12pt]{amsart}

 \usepackage{amsfonts,graphics,amsmath,amsthm,amsfonts,amscd,amssymb,amsmath,latexsym,multicol}
\usepackage{epsfig,url}
\usepackage{flafter}

\makeatletter

\def\jobis#1{FF\fi
  \def\predicate{#1}%
  \edef\predicate{\expandafter\strip@prefix\meaning\predicate}%
  \edef\job{\jobname}%
  \ifx\job\predicate
}

\makeatother

\if\jobis{proposal}%
\else
\fi

 \usepackage[matrix, arrow]{xy}

\DeclareMathOperator{\Supp}{Supp}

 \newcommand{\Q}{\mathbb Q}
 \newcommand{\R}{\mathbb R}
 
 \newcommand{\bir}{\dashrightarrow}

 \numberwithin{equation}{subsection}
 \numberwithin{footnote}{subsection}

 \newtheorem{cor}[subsection]{Corollary}

 \newtheorem{thm}[subsection]{Theorem}

{
 \newtheorem{defn}[subsection]{Definition}

 \newtheorem{rem}[subsection]{Remark}

}

 \newcommand{\ke}[1]{$\acute{\mbox{e}}$}
 \newcommand{\ku}[1]{$\acute{\mbox{u}$}}
 \newcommand{\kl}[1]{$\acute{\mbox{l}}$}
 \newcommand{\kh}[1]{$\acute{\mbox{h}}$}
 \newcommand{\kr}[1]{$\acute{\mbox{r}}$}
 \newcommand{\kx}[1]{$\acute{\mbox{x}}$}
 \newcommand{\ki}[1]{${\^\i}$}

\title{Log minimal models according to Shokurov}
\author{Caucher Birkar}
\date{\today}
\begin{document}
\maketitle

\begin{abstract}
Following Shokurov's ideas, we give a short proof of the following klt version of his result: termination of 
terminal log flips in dimension $d$ implies that any klt pair  
of dimension $d$ has a log minimal model or a Mori fibre space.  
Thus, in particular, any klt pair of dimension $4$ has a log minimal model or a Mori fibre space.
\end{abstract}


\section{Introduction}

All the varieties in this paper are assumed to be over an algebraically closed
field $k$ of characteristic zero. We refer the reader to section 2 for notation and terminology.

Shokurov [\ref{ordered}] proved that the log minimal model program (LMMP) in dimension $d-1$ and 
termination of terminal log flips in dimension $d$ imply existence of a log minimal model or a 
Mori fibre space for any lc pair of dimension $d$. Following Shokurov's method and using results of [\ref{BCHM}], we 
prove that termination of terminal log flips in dimension $d$ implies existence of a log minimal model or a 
Mori fibre space for any klt pair of dimension $d$. 

In this paper, by termination of terminal log flips in dimension $d$ we will mean termination of any sequence $X_i\bir X_{i+1}/Z_i$ of log flips$/Z$ starting with a $d$-dimensional klt pair $(X/Z,B)$ which is terminal in codimension $\ge 2$.

\begin{thm}\label{main}
Termination of terminal log flips in dimension $d$ implies that any klt pair $(X/Z,B)$ of dimension $d$ has a 
log minimal model or a Mori fibre space.
\end{thm}

\begin{cor}\label{c-dim4}
Any klt pair $(X/Z,B)$ of dimension $4$ has a log minimal model or a Mori fibre space.
\end{cor}

Note that, in the corollary, when $(X/Z,B)$ is effective (eg of nonnegative Kodaira dimension), log minimal models are constructed in [\ref{B2}] using different methods.

\section*{Acknowledgements}

I would like to thank V.V. Shokurov for useful comments.


\section{Basics}

Let $k$ be an algebraically closed field of characteristic zero. For an $\R$-divisor $D$ on a variety $X$ over $k$, we use $D^\sim$ to denote the birational transform of $D$ on a specified birational model of $X$. 

\begin{defn}
A pair $(X/Z,B)$ consists of normal quasi-projective varieties $X,Z$ over $k$, an $\R$-divisor $B$ on $X$ with
coefficients in $[0,1]$ such that $K_X+B$ is $\mathbb{R}$-Cartier, and a projective 
morphism $X\to Z$. $(X/Z,B)$ is called log smooth if $X$ is smooth and $\Supp B$ has simple normal crossing singularities. 

For a prime divisor $D$ on some birational model of $X$ with a
nonempty centre on $X$, $a(D,X,B)$
denotes the log discrepancy. $(X/Z,B)$ is terminal in codimension $\ge 2$ if $a(D,X,B)>1$ whenever $D$ is exceptional$/X$. Log flips preserve this condition 
but divisorial contractions may not.
\end{defn}

Let $(X/Z,B)$ be a klt pair. By a log flip$/Z$ we mean the flip of a $K_X+B$-negative extremal flipping contraction$/Z$. A sequence of log flips$/Z$ starting with $(X/Z,B)$ is a sequence $X_i\bir X_{i+1}/Z_i$ in which  $X_i\to Z_i \leftarrow X_{i+1}$ is a $K_{X_i}+B_i$-flip$/Z$ and $B_i$ is the birational transform of $B_1$ on $X_1$, and $(X_1/Z,B_1)=(X/Z,B)$. By termination of terminal log flips in dimension $d$ we mean termination of such a sequence in which $(X_1/Z,B_1)$ is a $d$-dimensional klt pair which is terminal in codimension $\ge 2$. Now assume that $G\ge 0$ is an $\R$-Cartier divisor on $X$. 
A sequence of $G$-flops$/Z$ starting with $(X/Z,B)$ is a sequence $X_i\bir X_{i+1}/Z_i$ in which  $X_i\to Z_i \leftarrow X_{i+1}$ is a $G_i$-flip$/Z$ such that $K_{X_i}+B_i\equiv 0/Z_i$ where $G_i$ is the birational transform of $G$ on $X=X_1$.\\

\begin{rem}\label{r-rays}
We borrow a result of Shokurov [\ref{ordered}, Corollary 9, Addendum 4] concerning extremal rays. 
Let $(X/Z,B)$ be a $\Q$-factorial klt pair and $F$ a reduced divisor on $X$. Then, there is $\epsilon>0$ such that 
if $G\ge 0$ is an $\R$-divisor supported in $F$ satisfying\\\\
(1) $||G||<\epsilon$ where $||.||$ denotes the maximum of coefficients, and\\
(2)  $(K_X+B+G)\cdot R<0$ for an extremal ray $R$,\\

then $(K_X+B)\cdot R\le 0$. This follows from certain numerical properties of log divisors such as [\ref{ordered}, Proposition 1] which is essentially the boundedness of the length of an extremal ray.    
Moreover, $\epsilon$ can be chosen such that for any $\R$-divisor $G'\ge 0$ supported in $F$ and any sequence $X_i\bir X_{i+1}/Z_i$ of 
$G'$-flops starting with $(X/Z,B)$ satisfying\\\\
(1') $||G_i||<\epsilon$ where $G_i\ge 0$ is a multiple of $G_i'$, the birational transform of $G'$, and\\
(2')  $(K_{X_i}+B_i+G_i)\cdot R<0$ for an extremal ray $R$ on $X_i$,\\

we have $(K_{X_i}+B_i)\cdot R\le 0$. 
In other words, $\epsilon$ is preserved after $G'$-flops but possibly only in the direction of $G'$. These claims are proved in [\ref{ordered}, Corollary 9, Addendum 4].
\end{rem}

\begin{defn}[Cf., {[\ref{B2}, \S 2]}]\label{d-mmodel}
Let $(X/Z,B)$ be a klt pair, $(Y/Z,B_Y)$ a $\Q$-factorial klt pair, $\phi\colon X\bir Y/Z$ a birational map such that $\phi^{-1}$ does not contract divisors, and $B_Y$ the birational transform of $B$. Moreover, assume that 
$$
a(D,X,B)\le a(D,Y,B_Y)
$$ 
for any prime divisor $D$ on birational models of $X$ and assume that the strict inequality holds for any prime divisor $D$ on $X$ which is exceptional/$Y$.

We say that $(Y/Z,B_Y)$ is a log minimal model of $(X/Z,B)$ if $K_Y+B_Y$ is nef$/Z$. 
On the other hand, we say that $(Y/Z,B_Y)$ is a Mori fibre space of $(X/Z,B)$ if 
there is a $K_Y+B_Y$-negative extremal contraction $Y\to Y'/Z$ such that $\dim Y'<\dim Y$. 
\end{defn}

Typically, one obtains a log minimal model or a Mori fibre space by a finite sequence of divisorial contractions and log flips. 

\begin{rem}\label{r-q-fact}
Let $(X/Z,B)$ be a klt pair and $W\to X$ a log resolution. Let $B_W=B^\sim+(1-\epsilon)\sum E_i$ where  $0<\epsilon\ll 1$ and $E_i$ are the exceptional$/X$ divisors on $W$. Remember that $B^\sim$ is the birational transform of $B$. If $(Y/X,B_Y)$ is a log minimal model of $(W/X,B_W)$, which exists by [\ref{BCHM}], then by the negativity lemma $Y\to X$ is a small $\Q$-factorialisation of $X$. To find a log minimal model or a Mori fibre space of $(X/Z,B)$, it is enough to find one for $(Y/Z,B_Y)$. So, one could assume that $X$ is $\Q$-factorial by replacing it with $Y$.
\end{rem}

Let $(X/Z,B+C)$ be a $\Q$-factorial klt pair such that $K_X+B+C$ is nef/$Z$. By [\ref{B2}, Lemma 2.6], either $K_X+B$ is nef/$Z$ or there is an extremal ray $R/Z$ such
that $(K_X+B)\cdot R<0$ and $(K_X+B+\lambda_1 C)\cdot R=0$ where
$$
\lambda_1:=\inf \{t\ge 0~|~K_X+B+tC~~\mbox{is nef/$Z$}\}
$$
and $K_X+B+\lambda_1 C$ is nef$/Z$. Now assume that $R$ defines a divisorial contraction or a log flip $X\bir X'$. We can consider $(X'/Z,B'+\lambda_1 C')$  where $B'+\lambda_1 C'$ is the birational transform of $B+\lambda_1 C$ and continue the argument. That is, either $K_{X'}+B'$ is nef/$Z$ or there is an extremal ray $R'/Z$ such
that $(K_{X'}+B')\cdot R'<0$ and $(K_{X'}+B'+\lambda_2 C')\cdot R'=0$ where
$$
\lambda_2:=\inf \{t\ge 0~|~K_{X'}+B'+tC'~~\mbox{is nef/$Z$}\}
$$
and $K_{X'}+B'+\lambda_2 C'$ is nef$/Z$. By continuing this process, we obtain a 
special kind of LMMP on $K_X+B$ which we refer to as the \emph{LMMP with scaling of $C$}. If it terminates, then we obviously get a log minimal model or a Mori fibre space for $(X/Z,B)$. Note that the 
required log flips exist by [\ref{BCHM}].


\section{Proofs}

\begin{proof}(of Theorem \ref{main})
Let $(X/Z,B)$ be a klt pair of dimension $d$. By Remark \ref{r-q-fact}, we can assume that $X$ is $\Q$-factorial. Let $H\ge 0$ be an $\R$-divisor which is big$/Z$ so that $K_X+B+H$ is klt and nef$/Z$. Run the LMMP/$Z$ on $K_X+B$ with scaling of 
$H$. If the LMMP terminates, then we get a log minimal model or a Mori fibre space. Suppose that 
we get an infinite sequence $X_i \bir X_{i+1}/Z_i$ of log flips$/Z$ where we may also assume that $(X_1/Z,B_1)=(X/Z,B)$. 

Let $\lambda_i$ be the threshold on $X_i$ determined by the LMMP with scaling as explained in section 2. So,   
$K_{X_i}+B_i+\lambda_i H_i$ is nef$/Z$, $(K_{X_i}+B_i)\cdot R_i<0$ and  $(K_{X_i}+B_i+\lambda_i H_i)\cdot R_i=0$ where $B_i$ and $H_i$ are the  birational transforms of 
$B$ and $H$ respectively and $R_i$ is the extremal ray which defines the flipping contraction $X_i\to Z_i$. Obviously, $\lambda_i\ge \lambda_{i+1}$.

Put $\lambda=\lim_{i\to \infty} \lambda_i$.  
If the limit is attained, that is, $\lambda=\lambda_i$ for some $i$, then the sequence terminates by [\ref{BCHM}, Corollary 1.4.2]. So, we 
assume that the limit is not attained. Actually, if $\lambda>0$, again [\ref{BCHM}] implies that the 
sequence terminates. However, we do not need to use [\ref{BCHM}] in this case. In fact, by replacing 
$B_i$ with $B_i+\lambda H_i$, we can assume that $\lambda=0$ hence $\lim_{i\to \infty} \lambda_i=0$. 

Put $\Lambda_i:=B_i+\lambda_i H_i$. Since we are assuming that terminal log flips terminate, or alternatively by [\ref{BCHM}, Corollary 1.4.3], we can construct a terminal (in codimension $\ge 2$) crepant model  
$(Y_i/Z,\Theta_i)$ of $(X_i/Z,\Lambda_i)$. A slight modification of the argument  in Remark \ref{r-q-fact} would do this. Note that we can assume that all the $Y_i$ are isomorphic to $Y_1$ in codimension 
one perhaps after truncating the sequence. Let $\Delta_1=\lim_{i\to \infty} \Theta_i^{\sim}$ on $Y_1$ and let $\Delta_i$ be its birational 
transform on $Y_i$. The limit is obtained component-wise. 

 Since $H_i$ is big$/Z$ and $K_{X_i}+\Lambda_i$ is klt and nef$/Z$, $K_{X_i}+\Lambda_i$ and $K_{Y_i}+\Theta_i$ are semi-ample$/Z$ by the base point freeness theorem for $\R$-divisors. 
Thus, $K_{Y_i}+\Delta_i$ is a limit of movable/$Z$ divisors which in particular means that it is pseudo-effective/$Z$. Note that if $K_{Y_i}+\Delta_i$ is not  pseudo-effective/$Z$, we get a contradiction by [\ref{BCHM}, Corollary 1.3.2].

Now run the LMMP$/Z$ on $K_{Y_1}+\Delta_1$. No divisor will be contracted again because  $K_{Y_1}+\Delta_1$ is a 
limit of movable/$Z$ divisors. Since $K_{Y_1}+\Delta_1$ is terminal in codimension $\ge 2$, by assumptions, the LMMP terminates with a log minimal model $(W/Z,\Delta)$. 
By construction, $\Delta$ on $W$ is the birational transform of $\Delta_1$ on $Y_1$ and $G_i:=\Theta_i^{\sim}-\Delta$ on $W$ satisfies $\lim_{i\to \infty} G_i=0$. 

By Remark \ref{r-rays}, for each $G_i$ with $i\gg 0$,  we can run the LMMP/$Z$ on $K_W+\Delta+G_i$ which will 
be a sequence of $G_i$-flops, that is, $K+\Delta$ would be numerically zero on all the extremal rays contracted in the 
process. No divisor will be contracted because $K_W+\Delta+G_i$ is movable$/Z$. 
The LMMP ends up with a log minimal model $(W_i/Z,\Omega_i)$. Here, $\Omega_i$ is the birational transform of $\Delta+G_i$ and so of $\Theta_i$. Let $S_i$ be the lc model of $(W_i/Z,\Omega_i)$ 
which is the same as the lc model of  $(Y_i/Z,\Theta_i)$ and that of $(X_i/Z,\Lambda_i)$ because $K_{W_i}+\Omega_i$ and $K_{Y_i}+\Theta_i$ are nef$/Z$ with  $W_i$ and $Y_i$ being isomorphic in codimension one, and $K_{Y_i}+\Theta_i$ is the pullback of $K_{X_i}+\Lambda_i$. Also note that since $K_{X_i}+B_i$ is pseudo-effective$/Z$, $K_{X_i}+\Lambda_i$ is big$/Z$ hence $S_i$ is birational to $X_i$.

By construction $K_{W_i}+\Delta^{\sim}$ is nef$/Z$ and it turns out that $K_{W_i}+\Delta^{\sim}\sim_{\R} 0/S_i$. Suppose that this is not the case. Then, $K_{W_i}+\Delta^{\sim}$ is not numerically zero$/S_i$ hence there is some 
curve $C/S_i$ such that $(K_{W_i}+\Delta^{\sim}+G_i^{\sim})\cdot C=0$ but $(K_{W_i}+\Delta^{\sim})\cdot C>0$ 
which implies that $G_i^{\sim}\cdot C<0$. Hence, there is a  $K_{W_i}+\Delta^{\sim}+(1+\tau)G_i^{\sim}$-negative extremal ray $R/S_i$ 
for any $\tau>0$. This contradicts Remark \ref{r-rays} because we must have 
$$
(K_{W_i}+\Delta^{\sim}+G_i^{\sim})\cdot R=(K_{W_i}+\Delta^{\sim})\cdot R=0
$$ 
Therefore, $K_{W_i}+\Delta^{\sim}\sim_{\R} 0/S_i$. Now $K_{X_i}+\Lambda_i\sim_{\R}0/Z_i$ implies that $Z_i$ is over $S_i$ and so  $K_{Y_i}+\Delta_i \sim_{\R} 0/S_i$. On the other hand, $K_{X_i}+B_i$ is the pushdown of $K_{Y_i}+\Delta_i$ hence $K_{X_i}+B_i\sim_{\R} 0/S_i$. Thus, $K_{X_i}+B_i\sim_{\R} 0/Z_i$ and this contradicts the fact that $X_i\to Z_i$ is a $K_{X_i}+B_i$-flipping contraction. So, the sequence of flips terminates and 
this completes the proof.
\end{proof}

\begin{proof}(of Corollary \ref{c-dim4})
Since terminal log flips terminate in dimension $4$ by [\ref{F}][\ref{Sh2}], the 
result follows from the Theorem.
\end{proof}


\vspace{2cm}

\flushleft{DPMMS}, Centre for Mathematical Sciences,\\
Cambridge University,\\
Wilberforce Road,\\
Cambridge, CB3 0WB,\\
UK\\
email: c.birkar@dpmms.cam.ac.uk

\end{document}